\documentclass[preprint,11pt]{elsarticle}

\usepackage{subfig,xcolor,graphicx}

\usepackage{fullpage}

\usepackage{amsmath,amssymb,amsthm}

\newtheorem{lemma}{Lemma}
\newtheorem{theorem}{Theorem}
\newtheorem{corollary}{Corollary}

\newcommand{\psiQex}{\psi_Q^{\rm explicit}}

\newcommand{\brev}{}

\newcommand{\Beginproof}{\vspace{0mm} \parindent=0pt
         {\bf Proof.} \hspace{3mm} \parindent=3ex}

\newcommand{\Endproof}{$\Box$ \vspace{5mm}
                        \parindent=3ex}

\newcommand{\Remark}{\vspace{0mm} \parindent=0pt
         {\bf Remark.} \hspace{0mm} \parindent=3ex}

\graphicspath{{}}

\begin{document}
\begin{frontmatter}
\title{Singular solutions of the subcritical nonlinear Schr\"odinger equation}
\author{Gadi Fibich}\ead{fibich@tau.ac.il}
\address{School of Mathematical Sciences, Tel Aviv University, Tel Aviv 69978, Israel}
\begin{abstract}
{\brev
We show that the subcritical $d$-dimensional nonlinear Schr\"odinger equation  $i \psi_t + \Delta \psi + |\psi|^{2 \sigma} \psi = 0$, where
$1<\sigma d<2$, admits smooth solutions that become singular in~$L^p$ for $p^*<p \le \infty$, where $p^*:=\frac{\sigma d}{\sigma d -1}$.
Since $\lim_{\sigma d \to 2-} p^* = 2$, these solutions can collapse at any $2<p \le \infty$, and in particular for $p = 2 \sigma+2$.  
}
\end{abstract}
\end{frontmatter}

\section{\label{sec:intro}Introduction}

The focusing nonlinear Schr\"odinger equation (NLS)
\begin{equation}
   \label{eq:NLS_sc_explicit}
   i \psi_t(t, {\bf x}) + \Delta \psi + |\psi|^{2 \sigma} \psi = 0,  \qquad \psi(0, {\bf x})  = \psi_0({\bf x}),
\end{equation}
where~$	{\bf x} =\left(x_1,\dots,x_d\right)\in {\Bbb R}^d$ and~$\Delta = \sum_{j=1}^d\partial_{x_j x_j}$ is the Laplacian, 
has been the subject of intense study, due to its role in
various areas of physics, such as nonlinear optics and Bose-Einstein Condensates
(BEC).  The NLS is called subcritical, critical, and supercritical if $\sigma d<2$, $\sigma d=2$, and $\sigma d>2$, respectively.  It is well-known that in the critical and supercritical cases, the NLS~\eqref{eq:NLS_sc_explicit} possesses solutions that become singular in a finite time~\cite{Sulem-99}. In this study we show that, contrary to common belief, the 
subcritical NLS also admits solutions that become singular in a finite time.

 Most of NLS theory has been developed for solutions that are in~$H^1({\Bbb R}^d)$. In this case, the initial condition $\psi_0 \in H^1$, and the NLS solution is said to become singular at $t = T_c$, if $\psi(t) \in H^1$ for $0 \le t< T_c$, and 
$\lim_{t \to T_c} ||\psi(t)||_{H^1} = \infty$.
 In 1983, Weinstein proved that all $H^1$~solutions of the subcritical NLS exist globally: 
\begin{theorem}[\cite{Weinstein-83}]
   \label{thm:existence}
  Let $\psi$ be a solution of the NLS~\eqref{eq:NLS_sc_explicit}, let $0<\sigma d<2$, and let $\psi_0 \in H^1$.
Then, $\psi$~exists globally in~$H^1$.
\end{theorem} 
  Until now, this result has been interpreted as implying that the subcritical NLS does not admit singular solutions.
  In this study we show that if we do not restrict ourselves to $H^1$~solutions, then the subcritical NLS also admits singular solutions.
 Here, by singular we mean that {\brev there exists some $2<p<\infty$, such that~$||\psi||_p$ becomes infinite in a finite time}.\footnote{In the case of 
$H^1$~solutions of the NLS, blowup of the $H^1$ norm implies  blowup of {\brev the $L^p$ norms for $\brev 2 \sigma+2 \le p \le \infty$, see~\ref{app:blowup-Linfty}.}}
Our main result is as follows:
\begin{theorem}  
   \label{thm:collapse-subcritical-main}
{\brev Let }
\begin{equation}
\label{eq:p-range}
{\brev p^* <p \le \infty, \qquad p^*:=\frac{\sigma d}{\sigma d -1}.}
\end{equation}
Then, the subcritical NLS with $1<\sigma d<2$ admits classical solutions that becomes singular 
at a finite time~$T_c$ {\brev in~$L^p$}, i.e., 
$$
||\psi(t)||_p< \infty, \qquad 0 \le t < T_c,
$$
and 
$$
\lim_{t \to T_c}  ||\psi (t)||_p  =  \infty.
$$ 
\end{theorem}

Theorem~\ref{thm:collapse-subcritical-main} follows from the following Theorem:

\begin{theorem}  
   \label{thm:collapse-subcritical}
    {\brev Let $p$ be in the range~\eqref{eq:p-range}}, let $1<\sigma d<2$, 
let $a>0$ be a positive constant, and let  
 $Q(\rho)$ be the solution of 
\begin{subequations}
\label{eq:Q-Zakharov}
\begin{equation}
 \Delta Q(\rho) - Q + ia \left(\frac{1}{\sigma} Q+ \rho Q' \right) + |Q|^{ 2 \sigma} Q = 0, 
   \qquad 0< \rho<\infty, 
\end{equation}
subject to 
\begin{equation}
\label{eq:Q-Zakharov-ic}
  0 \not= Q(0) \in {\Bbb C}, \qquad Q'(0) = 0.  
\end{equation}
\end{subequations}
Let 
\begin{subequations}
   \label{eq:psiQexplicit-ABC}
\begin{equation}
   \label{eq:psiQexplicit}
  \psiQex(t,r) =  \frac{1}{L^{1/\sigma}(t)} Q(\rho) e^{i \tau(t)},
\end{equation}
where 
\begin{equation}
  \label{eq:psiQexp-Zakharov-L}
r = |{\bf x}|, \qquad    L(t) = \sqrt{2a (T_c-t)}, 
\end{equation}
and
\begin{equation}
  \label{eq:psiQexp-Zakharov-rho}
  \rho  = \frac{r}{L(t)}, \qquad \tau  = \int_0^t \frac{1}{L^2(s)} \, ds = \frac{1}{2 a} \log \frac{T_c}{T_c-t}.
\end{equation}
\end{subequations}
Then, $\psiQex$ is an explicit solution of the 
 subcritical NLS that becomes singular in~$\brev L^p$
as $ t \longrightarrow T_c$. 
\end{theorem}

\Remark Although $Q$, hence also $\psiQex$, is not in~$H^1$, it is smooth, and it decays to zero as $|{\bf x}|  \longrightarrow \infty$, see Lemma~\ref{lem:Q-subcritical}.

{\brev Since 
$$\lim_{\sigma d \to 2-}p^* = 2+,
$$
then for any $2<p<\infty$, there exists a singular solution of a subcritical NLS that becomes singular in $L^p$. In particular, if $\sigma d$ is sufficiently close to~2 from below, then $\psiQex$ becomes singular in $L^{2 \sigma +2}$.    

\Remark  The linear Schr\"odinger equation   $i \psi_t + \Delta \psi  = 0$ admits the fundamental solution 
$ \psi = \frac{1}{(4 \pi i t )^{d/2}} e^{i|{\bf x}|^2/4 t }$,
which becomes singular in finite time in~$L^\infty$~\cite{Cordero-08}. 
Unlike~$\psiQex$, however, this solution does not become singular in~$L^p$ for any finite~$p$.

}

\section{Proof of Theorem~\ref{thm:collapse-subcritical}}

We begin with the following result.
\begin{lemma}
  \label{lem:psiQexp-Zakharov}
let $\psiQex$ be defined as in Theorem~\ref{thm:collapse-subcritical}.
Then, $\psiQex$ is an explicit solution of the NLS~\eqref{eq:NLS_sc_explicit}.
\end{lemma}
\Beginproof Substituting~$\psiQex$ in the NLS~\eqref{eq:NLS_sc_explicit} and carrying out the differentiation 
proves the result. \Endproof

The result of Lemma~\ref{lem:psiQexp-Zakharov} was used by Zakharov~\cite{Zakharov-84}, and subsequently by others (see~\cite{Sulem-99} and the references therein), in the study of singular $H^1$~solutions of the supercritical NLS. These solutions
undergo a quasi self-similar collapse, in which $\psiQex$ is the {\em asymptotic} blowup profile of the collapsing core of the solution. 
Here, in contrast, $ \psiQex$ is an explicit,  ``truly'' self-similar solution of the subcritical NLS.

We now establish the decay at infinity of all solutions of equation~\eqref{eq:Q-Zakharov}:
\begin{lemma}
   \label{lem:Q-subcritical}
Let $a>0$ and let $1<\sigma d<2$. Then, for any $Q(0) \in {\Bbb C}$, 
the solution of equation~\eqref{eq:Q-Zakharov} exists, is unique,  and decays to zero as $\rho \longrightarrow \infty$,
so that
$$
  |Q| = O(\rho^{-d+1/\sigma}), \qquad \rho \longrightarrow \infty. 
$$ 
{\brev Therefore, $Q$ is in~$L^p$ for any $\brev p^*<p \le \infty$.}
\end{lemma}
\Beginproof  The proof is nearly identical to the proof of Johnson and Pan in the supercritical case~\cite{Johnson-93}, see~\ref{app:Q-subcritical}.
\Endproof


\begin{lemma}
{\brev For any $\brev p^*<p \le \infty$}, 
  $\psiQex$ becomes singular in~$L^{p}$
as $\left. t \longrightarrow T_c \right.$. 
\end{lemma}
\Beginproof
Since 
$$
 ||\psiQex (t)||_p  =  
\frac{||Q ||_p}{L^{1/\sigma}(t)},
$$
and $||Q||_p<\infty$, the result follows.
\Endproof

This concludes the proof of Theorem~\ref{thm:collapse-subcritical}.

\section{The $Q$ equation in the subcritical case}

As in~\cite{LeMesurier-88a}, 
the far-field asymptotics of~$Q$ can be calculated using the WKB method:
\begin{lemma}
   \label{lem:Q-rho>>1-subcritical}
Let $Q(\rho)$ be a solution of equation~\eqref{eq:Q-Zakharov}, where 
$1<\sigma d <2$.
Then,
\begin{equation}
   \label{eq:Q=c1Q1+c2Q2}
   Q(\rho) \sim c_1 Q_1(\rho) + c_2 Q_2(\rho), \qquad  \rho \longrightarrow \infty, 
\end{equation}
where $c_1$ and $c_2$ are complex numbers, and
$$
 Q_1  \sim  \rho^{-i/a-1/\sigma} , \qquad  Q_2  \sim  e^{-i a \rho^2/2}   \rho^{i/a-d+1/\sigma} .
$$
\end{lemma}
\Beginproof See~\ref{app:Q-rho>>1-subcritical}.
\Endproof

\begin{corollary}
   \label{cor:Q1Q2in}
If $1<\sigma d <2$, then  
$$
 Q_1 \in L^2({\Bbb R}^d), \qquad \nabla Q_1 \in L^2({\Bbb R}^d), 
$$
and                
$$
 Q_2 \not\in L^2({\Bbb R}^d), \qquad \nabla Q_2 \not\in L^2({\Bbb R}^d).
$$
In addition,  $Q_1 \in L^{p}({\Bbb R}^d)$ for any $2 \le p \le \infty$, and  
$$
 Q_2 \in L^{p}({\Bbb R}^d), \qquad \frac{\sigma d}{\sigma d -1} <p \le \infty.
$$
\end{corollary}
\Beginproof This follows from Lemma~\ref{lem:Q-rho>>1-subcritical}.
\Endproof


%



In the supercritical case, a key role is played by the zero-Hamiltonian solutions of the $Q$~equation, which 
behave as~$c_1 Q_1$ at large~$\rho$~\cite{Sulem-99}. We now show that there are no such solutions in the subcritical case: 
\begin{lemma}
   \label{lem:no-H1-Qsols-suncritical}
 When $1<\sigma d <2$, there are no nontrivial solutions of the $Q$~equation~\eqref{eq:Q-Zakharov}, such that $c_2 = 0$, i.e., that
$$
   Q(\rho) \sim c_1 Q_1(\rho),  \qquad  \rho \longrightarrow \infty.
$$
\end{lemma}
\Beginproof
  By negation. Assume that there is such a~$Q$. In this case, it follows from Corollary~\ref{cor:Q1Q2in} that  
$Q \in H^1$. Hence, $\psiQex$ is a solution of the subcritical NLS that becomes singular 
in~$H^1$, which is in contradiction with Theorem~\ref{thm:existence}.~\Endproof

\section{Simulations}

\begin{figure}[ht]
    \centerline{\scalebox{0.8}{\includegraphics{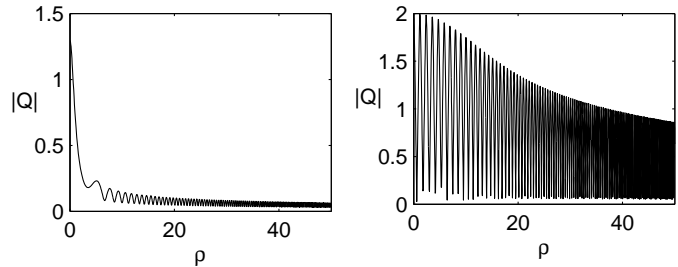}}}
 \caption{Solutions of equation~\eqref{eq:Q-Zakharov} with $d=1$, $\sigma = 1.9$,  $a=0.5145$, and with A:~$Q(0) = 1.2953$, 
and B:~$Q(0) = 3$.  
   }
   \label{fig:Q-subcritical}
\end{figure}

  Figure~\ref{fig:Q-subcritical} shows two numerical solutions of equation~\eqref{eq:Q-Zakharov}.
As expected, see Lemma~\ref{lem:Q-rho>>1-subcritical},
$$|Q| \sim 
   \big|c_1  \rho^{-i/a-1/\sigma} + c_2 e^{-i a \rho^2/2}   \rho^{i/a-d+1/\sigma}\big|
$$ 
decreases to zero  as $\rho \longrightarrow \infty$, while undergoing faster and faster oscillations. 
The ``cleaner picture'' in Figure~\ref{fig:Q-subcritical}A has to do with the fact that
the values of~$a$ and~$Q(0)$ were chosen so as to minimize the value of~$c_2$.\footnote{These values of~$a$ and~$Q(0)$ were calculated using the shooting algorithm of Budd, Chen and Russel~\cite{Budd-99} for calculating the zero-Hamiltonian solutions in the supercritical case. }\footnote{The value of~$c_2$ cannot be equal to zero, see Lemma~\ref{lem:no-H1-Qsols-suncritical}.}

\begin{figure}[ht]
    \centerline{\scalebox{0.7}{\includegraphics{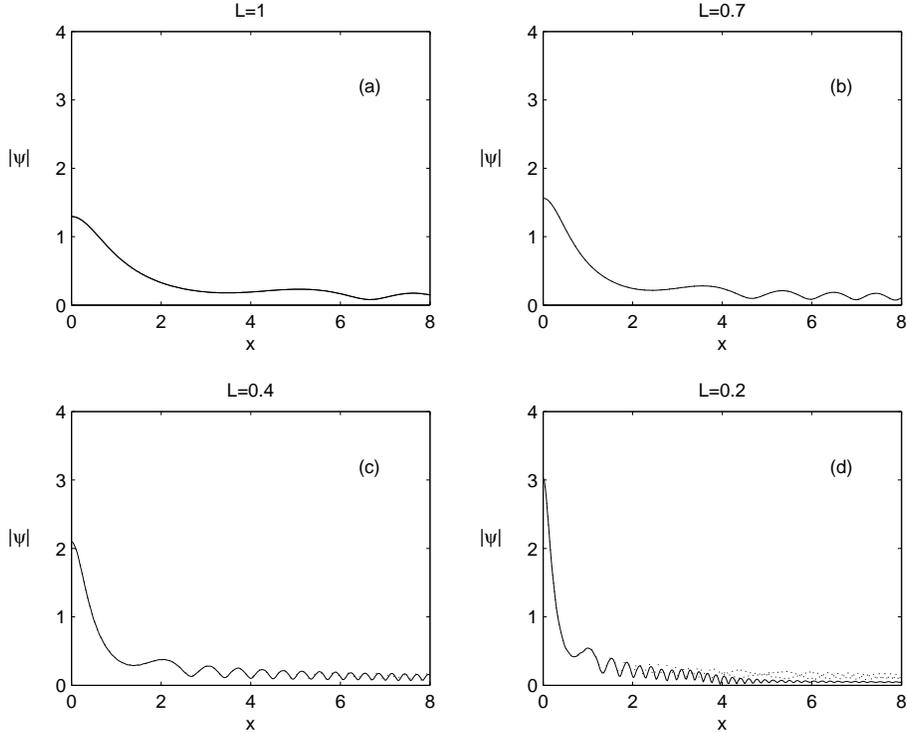}}}
 \caption{Numerical solution of the subcritical NLS with the initial condition equation~$\psi_0(x) =Q(x)$ (solid). Dotted line is the analytic solution~$\psiQex$. a)~$t=0$, $L=1$; b)~$t=0.4956$, $L=0.7$; c)~$t= 0.8163$, $L=0.4$;  
d)~$t= 0.9329$, $L=0.2$.
   }
   \label{fig:QprofileSelfSimilar}
\end{figure}

In Figure~\ref{fig:QprofileSelfSimilar} we solve numerically the subcritical NLS with $d=1$, $\sigma = 1.9$, and the initial condition
$
  \psi_0(x) =Q(x),
$   
where $Q$ is taken from Figure~\ref{fig:Q-subcritical}A.  By Lemma~\ref{lem:psiQexp-Zakharov}, the analytic solution of this equation is given by~$\psiQex$ with $L = \sqrt{1-2 a t}$. As expected, the numerical solution agrees with the analytic solution, thus providing the first ever simulation of a singular solution of the subcritical NLS. 

{\brev In these simulations, we used a standard fourth-order finite-difference implicit scheme with~$dx = 0.05$ and~$dt = 0.001$
  over the spatial domain  
$-70 \le x \le 70$. Nevertheless, because of the slow decay and the ever faster oscillations as $x \longrightarrow \infty$, the 
agreement between the analytic and numerical solutions breaks down after focusing by less than~3 (see Figure~\ref{fig:QprofileSelfSimilar}d). 
 We could, of course, take an even larger domain with a finer mesh.
In that case, the numerical solution would simply bifurcate from the analytic one at a higher focusing level. The point of this simulation, however, is not to establish numerically the existence of a singular subcritical solution (which we prove rigorously), but rather to illustrate the numerical difficulties in computing this solution, by showing that even with a relatively large domain and a fine grid, the numerical solution breaks down after focusing by less than~3. This suggests that numerical simulations may be useless in studying the stability of these solutions.  
}
%

\section{Final remarks}

    Until now, it was believed that only the critical and supercritical NLS admits singular solutions.
      In this study we showed that if we do not limit ourselves to $H^1$~solutions, then the subcritical NLS also admits solutions that become singular  at a finite time. 
This finding raises several questions, which are currently open. One question is whether the explicit singular solutions are stable. As note, 
this question is hard to study numerically, because of the slow decay, coupled with the ever faster oscillations, of the solution at infinity. 
Another open question is whether the subcritical NLS admits singular solutions that are not self-similar.
The answers to these questions will determine whether singularity formation in the subcritical NLS will remain as an anecdote, 
or lead to  a new line of research.

%
%
%

\subsubsection*{Acknowledgments.}
This research was partially supported by the Israel Science
Foundation (ISF grant No. 123/08). We would like to thank Moran Klein for the NLS simulations of Figure~\ref{fig:QprofileSelfSimilar}.

\bibliographystyle{unsrt}
\bibliography{NLS}

\appendix

\section{Blowup of $\brev ||\psi||_{p}$}
  \label{app:blowup-Linfty}

The NLS conserves the power (mass, $L^2$ norm) and the Hamiltonian, i.e.,
$$
||\psi||_2^2 \equiv ||\psi_0||_2^2, \qquad 
  H(t) := ||\nabla \psi||_2^2 - \frac{1}{\sigma+1} ||\psi||_{2 \sigma+2}^{2 \sigma+2}  \equiv H(0).
$$  
Therefore, when  $||\psi||_{H^1}$ becomes infinite, then so does $||\nabla \psi||_{2}$, hence
$||\psi||_{2 \sigma+2}$. Therefore, since $||\psi||_2$ is conserved, it follows from the interpolation inequality for $L^p$~norms 
that {\brev  $ \brev ||\psi||_p$ also becomes infinite for $\brev 2 \sigma+2 \le p\le  \infty$.}

\section{Proof of Lemma~\ref{lem:Q-subcritical}}
  \label{app:Q-subcritical}

As in the proof of Johnson and Pan in the supercritical case~\cite{Johnson-93},
let 
$$
   Q(\rho) = u(\rho) e^{-i a \rho^2/4}, \qquad 
   u(\rho) = u_1(\rho)+ i u_2(\rho),
$$
where $u_1$ and~$u_2$ are real functions, let 
$$
  v_j(\rho) = \rho^{(d-1)/2} u_j(\rho), \qquad j = 1,2,
$$ 
let
$$
 t = \rho^2, \qquad x_j(t) = v_j(\rho), \qquad y_j(t) = \frac{dx_j}{dt},
$$
let
$$
 f_j(t) = t^{1/4} x_j(t), \qquad 
 g_j(t) = \frac{df_j}{dt}, 
$$ 
and let
$$
   H(t) = \frac{1}{2} (g_1^2+ g_2^2) + \frac{1}{8} \left( \lambda -\frac{1}{t} -\frac{e}{t^2} \right)   (f_1^2+ f_2^2)
   + \frac{1}{4(2 \sigma+2)} t^{-\beta} (f_1^2+ f_2^2)^{\sigma+1},
$$
where
$$
  \lambda = \frac{a^2}{4}, \qquad 
 \beta = 1+ \frac{\sigma d}{2}, \qquad e = \frac{1}{4} d (d-4).
$$
Then,
$$
  H'(t) = \frac{B}{4 t} (f_1 g_2-f_2 g_1) + \frac{1}{8 t^2}  \left( 1+ \frac{2e}{t} \right)  (f_1^2+ f_2^2)
-\frac{\beta}{4(2 \sigma+2)} t^{-\beta -1}  (f_1^2+ f_2^2)^{\sigma+1},
$$
where
$$
   B = a \left( \frac{d}{2}-\frac{1}{\sigma}  \right)<0.
$$
In addition, from the Cauchy-Schwartz inequality~\cite{Johnson-93},
$$
  |f_1 g_2| \le \frac{1}{2} \left[ \frac{\sqrt{\lambda}}{2} f_1^2 +  \frac{2}{\sqrt{\lambda}} g_2^2 \right],
\qquad 
  |f_2 g_1| \le \frac{1}{2} \left[ \frac{\sqrt{\lambda}}{2} f_2^2 +  \frac{2}{\sqrt{\lambda}} g_1^2 \right].
$$ 

Since $\beta>0$ and $B<0$,
\begin{eqnarray*}
  H'(t) &\le&  \frac{|B|}{4 t} (|f_1 g_2| + |f_2 g_1|) + \frac{1}{8 t^2}  \left( 1+ \frac{2e}{t} \right)  (f_1^2+ f_2^2).
 \\ 
 &\le &
\frac{|B|}{8 t} \left(\frac{\sqrt{\lambda}}{2} (f_1^2+f_2^2) + 
 \frac{2}{\sqrt{\lambda}} (g_1^2+g_2^2) \right) + \frac{1}{8 t^2}  \left( 1+ \frac{2e}{t} \right)  (f_1^2+ f_2^2)
 \\ 
 &= &
\frac{|B|}{2 \sqrt{\lambda} t} \left(\frac{\lambda}{8} (f_1^2+f_2^2) + 
 \frac{1}{2} (g_1^2+g_2^2) \right) + \frac{1}{8 t^2}  \left( 1+ \frac{2e}{t} \right)  (f_1^2+ f_2^2)
 \\ 
 & \le  &
\frac{|B|}{2 \sqrt{\lambda} t} \left( 1 + O\left(\frac{1}{t}\right) \right) H(t).
\end{eqnarray*}
Since $H(t)>0$ for large~$t$,
$$
   \frac{H'}{H} \le \frac{|B|}{2 \sqrt{\lambda} t}  + O\left(\frac{1}{t^2}\right).
$$
Therefore, as in~\cite{Johnson-93}, 
there exists a constant $c>0$, such that
$$
   H(t) \le c (1+ t^{2 |\alpha|}), \qquad 0 \le t<\infty,
$$
where
$$
  2 |\alpha| =\frac{|B|}{2 \sqrt{\lambda}} = \frac{1}{\sigma}-\frac{d}{2}. 
$$
Hence,
$$
   |f_j(t)|  \le   c (1+ t^{|\alpha|}),
$$
$$
   |x_j(t)|  \le   c t^{-1/4} (1+ t^{|\alpha|}), 
$$
$$
|v_j(\rho)| \le  c \rho^{-1/2} (1+ \rho^{2|\alpha|}),
$$
and
$$
 |u_j(\rho)|  \le c \rho^{-d/2} (1+ \rho^{2|\alpha|}).
$$
Therefore, 
$$
  |u| = O(\rho^{-d+1/\sigma}), \qquad \rho \longrightarrow \infty. 
$$

\Remark The only difference from the original proof of Johnson and Pan is that in in the supercritical case~$B>0$.
Therefore, we take the absolute value of~$B$, instead of~$B$, in the bounds for $H'$.

\section{Proof of Lemma~\ref{lem:Q-rho>>1-subcritical}}
   \label{app:Q-rho>>1-subcritical}
Let
$$
   Q(\rho)  =  e^{ -\frac{1}{2} \int \left( \frac{d-1}{\rho} +  ia \rho\right) }  Z(\rho)
 =  e^{ -i a \rho^2/4}  \rho^{-(d-1)/2}  Z(\rho).
$$
Therefore, the equation for~$Z$ is given by 
\begin{equation}
   \label{eq:Z-subcritical}
  Z''(\rho) + \left(\frac{a^2}{4} \rho^2 -1 -  i a \frac{d \sigma-2}{2 \sigma}  - \frac{(d-1)(d-3)}{4 \rho^2}
+ \left|Q \right|^{2 \sigma} \right) Z  = 0.
\end{equation}
Since by Lemma~\ref{lem:Q-subcritical}, $\lim_{\rho \to \infty} Q  = 0$,
%
 let us look for an asymptotic solution of the form 
$$
  Z = e^{w(\rho)}, \qquad w(\rho) \sim w_0(\rho)+ w_1(\rho) + \cdots.
$$
The equation for $\{w_i(t)\}$ is given by
\begin{equation}
 \label{eq:wi-sc}
(w_0''+ w_1'' + \cdots) + (w_0'+ w_1' + \cdots)^2 
+ \left(\frac{a^2}{4} \rho^2 -1 -  i a \frac{d \sigma-2}{2 \sigma}  - \frac{(d-1)(d-3)}{4 \rho^2} + \left|Q \right|^{2 \sigma} \right) = 0. 
\end{equation}
A-priori, the equation for the leading-order terms is 
$$
w_0'' + (w_0')^2 + \frac{a^2}{4} \rho^2 = 0.
$$
 The substitution $w_0 = c \rho^n$ shows that the order of the terms in this equation
is $\rho^{n-2}$, $\rho^{2n-2}$, and $\rho^2$, respectively. 
Since the only consistent way to balance the leading-order terms is if $n=2$,
the equation for the leading-order terms is given by 
$$
(w_0')^2 + \frac{a^2}{4} \rho^2 = 0.
$$
Therefore,
$$
 w_0' = \pm \frac{i a}{2} \rho, \qquad 
 w_0 = \pm \frac{i a}{4} \rho^2.
$$
The balance of the next-order terms is given by 
$$
 w_0''+ 2 w_0' w_1' - 1 -  i a \frac{d \sigma-2}{2 \sigma}  = 0.
$$
Substituting $w_0' = \pm i a  \rho/2$ and rearranging gives,
$$
  w_1' =  \mp\frac{ i}{a \rho}  \pm \frac{d \sigma-2}{2 \sigma} \frac{1}{\rho}  - \frac{1}{ 2 \rho},  
\qquad 
  w_1 =  \left( \mp\frac{ i}{a}  \pm \frac{d \sigma-2}{2 \sigma}  - \frac{1}{ 2} \right) \log \rho.  
$$
We will now show that $w_2 = o(1)$. Therefore, we obtained the two solutions 
\begin{eqnarray*}
 w^{(1)}(\rho) & = & i a \frac{\rho^2}{4} +  \left( -\frac{i}{a}  -\frac{1-d}{2} - \frac{1}{\sigma} \right)\log \rho + o \left(1 \right), 
\\
 w^{(2)}(\rho) & = & -i a \frac{\rho^2}{4} +\left( \frac{i}{a} + \frac{-1-d}{2} + \frac{1}{\sigma} \right)\log \rho + o \left(1 \right).
\end{eqnarray*}
Substituting $Q_i(\rho)  =  e^{ -i a \rho^2/4}  \rho^{-(d-1)/2}   e^{w^{(i)}(\rho)}$  leads to the result.

   In order to confirm that $w_2 = o(1)$, we note that the equation for~$w_2$ is given by
\begin{equation}
 \label{eq:w2-sc}
   w_1''+ (w_1')^2+ 2 w_0' w_2' - \frac{(d-1)(d-3)}{4 \rho^2}  + |Q|^{2 \sigma} = 0.
\end{equation}
In the case of $Q_1$, since $|Q_1|^{2  \sigma} \sim \rho^{-2}$, 
substituting the expressions for $w_0$ and $w_1$ gives
$$
  w_2' = O \left(\frac{1}{\rho^3} \right), \qquad  w_2 = O \left(\frac{1}{\rho^2} \right).
$$
In the case of $Q_2$, since 
$|Q_2|^{2  \sigma} \sim \rho^{-2 \sigma d +2} \gg \rho^{-2}$, 
 the leading-order equation for~$w_2$  becomes 
$$
   2 w_0' w_2' + |Q_2|^{2 \sigma}  = 0. 
$$
   Since $w_0'  \sim \rho$, then
$w_2' \sim  \rho^{-2 \sigma d +1}$ and 
$w_2 \sim  \rho^{-2 \sigma d +2} = o(1)$.

Finally, we note that this this proof is rigorous, since solutions of linear ODEs always have their asymptotics obtained by WKB calculations,
and the ODE~\eqref{eq:Z-subcritical} for~$Z$ is ``linear'', since it can be written as 
$$
  Z''(\rho) + \left(\frac{a^2}{4} \rho^2 -1 -  i a \frac{d \sigma-2}{2 \sigma}  - \frac{(d-1)(d-3)}{4 \rho^2}
+ O(\rho^{-d+1/\sigma}) \right) Z  = 0.
$$
\end{document}